# Mathematical analysis in high school : a fundamental dilemma

Carl Winsløw, University of Copenhagen

**Introduction**

Analysis, as construed here, is a domain of mathematics which treats problems related to limits, real and complex functions, and linear operators. While some of these problems have been known for thousands of years, the fundamentals of contemporary analysis – which include a rigorous theory of real numbers – have been established over the past 400 years. Analysis is closely linked to geometry and algebra, and also to a number of domains in the natural and social sciences. In particular, theoretical constructs like derivative and integral are historically linked to fundamental notions in mechanics and geometry (such as speed and area) while, today, derivatives and integrals are used in many other contexts.

The introduction into secondary level mathematics of elementary analysis, especially differential and integral calculus, was historically justified by the manifest and increasing importance of these elements both in pure mathematics and in other disciplines. How and when it was done clearly varies from one national or regional context to another, for instance it remains an option in the USA (cf. e.g. Spresser, 1979). In Denmark, the first timid introduction of "infinitely small and large quantities" as a mandatory topic in high school came as early as 1906. The teaching of "infinitesimal calculus" – that is, differential and integral calculus – became mandatory in the scientific stream from 1935. Ever since, the investigation of functions based on derivatives and integrals has remained a relatively stable and central part of the tasks posed for the national written examinations in the scientific streams of Danish high school (Petersen and Vagner, 2003). And it certainly remains a central element of the more advanced mathematics curriculum at this level. Most of the analysis exercises from the national final exams of the late 1930's could be found in today's exams, except for details of formulation.

Despite the stability of the "core" types of tasks – such as determining the extreme values of a given function on an interval – one may nevertheless point out two major periods of change, which are not specific to Danish high school but can be found in many other European countries:

- Around 1960, the range and formality of mathematical themes was significantly extended, especially in adjacent domains such as set theory and algebra – but all or most of the extensions were subsequently abandoned after a decade or two ;
- From around 1980, the progressive introduction of calculating devices in secondary schools has increasingly affected the teaching of certain core techniques in analysis.

In this chapter, we will first provide a theoretical framework for analysing and comparing different forms of organizing introductions to mathematical analysis, then illustrate it by two characteristic examples from the above periods of change as they occurred in Denmark, based on the national exam tasks and text books used in the two periods. We conclude by extracting from this a fundamental dilemma for the teaching of analysis at secondary level in view of the requirements and affordances provided by computer algebra systems on the one hand, and contemporary utilitarian school pedagogies on the other hand.

**An epistemological reference model**

As affirmed already in the first phrase of this chapter, the notion of limit is fundamental to mathematical analysis and in particular to elements that appear in most introductory calculus teaching. Among these elements, the most important is probably the notion of derivative function, and in order to present an explicit an general definition of derivation, any calculus text will have to introduce at least an informal explanation of what $\lim_{x \to a} f(x)$ means for a function $f$ defined in a neighborhood of $a$ (except possibly at $a$). Of course, the actual computation of such limits may also be useful for investigations of the function $f$ itself, as well as in other contexts. As a result, most calculus texts and syllabi include at least a little practice and theory related to limits, prior to the introduction of derivatives.

For their study of the teaching of limits of functions in Spanish high school, Barbé, Bosch, Espinoza and Gascón (2005) proposed an epistemological reference model to trace the didactic transposition of pertinent knowledge whose end result is the didactic process observed in the classroom. As this paper has appeared in a widely accessible journal, we recall only one main point while we use the full theoretical framework explained in that paper, in particular the basic notions of the anthropological theory of the didactical explained in section 2 of the paper.

The main point we wish to emphasise is that the authors identify two local mathematical organisations which they use for their study of the different stages of the didactic transposition of the basic theory of limits in a Spanish high school class. These organisations are:

- $MO_1$, termed "the algebra of limits", where the practice is unified by a discourse about how to compute limits in a variety of cases, and the practice blocks amount to such cases: each consists of a type of task with a technique that allows to solve all tasks of the given type – for instance, to compute $\lim_{x \to a} f(x)$ when $f$ is a polynomial, the answer is simply $f(a)$. At the theoretical level of this organisation there are algebraic rules like $\lim_{x \to a}(f(x) + g(x)) = \lim_{x \to a} f(x) + \lim_{x \to a} g(x)$ which are not further justified. Also the existence of limits is not problematised beyond the possibility of computation.
- $MO_2$, termed "the topology of limits", where the practice is unified by an abstract discourse and theory about limits, including for example a rigorous definition of what it means for $\lim_{x \to a} f(x)$ to exist; the types of tasks in this organisation include determining if a given function is differentiable at a given point, and to justify calculation rules like $\lim_{x \to a}(f(x) + g(x)) = \lim_{x \to a} f(x) + \lim_{x \to a} g(x)$ under appropriate assumptions.

The link between the two organisations is, at least in principle, clear: the practice block of $MO_2$ is needed to justify the theoretical level of $MO_1$ in a wider theoretical context (while, locally, the calculations rules for limits might be regarded as a kind of self-evident axioms). In this wider theoretical sense of limits, namely that of academic mathematics, one might even say that $MO_2$ comes first: before calculating $\lim_{x \to a} f(x)$ we need to define what it means, and that certainly includes non-trivial conditions for existence.

In didactic practice $MO_2$ does not need to come first. It is apparent from the study cited above, as well as from other research on the teaching of limits (often with less explicit reference models) that the practice block of $MO_1$ may in fact be taught and learned with relative ease and efficiency, with a theoretical block that is limited to informal and practice oriented explanation of the calculation rules. On the other hand, the teaching of $MO_2$ is usually absent or sparse at the secondary level both in Spain and elsewhwere. In fact, convergence is often described informally, based on examples of function graphs and verbal explanations of how the function value gets close to a limit value as the free variable "moves" towards a given value. The development of a practice block (with mathematical techniques related to $MO_2$) for students is quite rare; it would, for instance, imply giving students rigorous techniques to decide on the question of existence of $\lim_{x \to a} f(x)$ in concrete and non-trivial cases. By rigorous, we mean that the technique can be explained and justified at the theory level of $MO_2$, such as the example given in Barbé et al. (2005, p. 243). When it is done, it is often prepared by introducing first the simpler theory of *limits of sequences* of real numbers.

We note here, in passing, the strange and almost circular use of the term continuity found in the Spanish high school (see ibid., p. 255) and most likely in many similar institutions. The meaningfulness of this notion seems to be particularly affected by the lack of a practical block in the didactic transposition of $MO_2$, both in the prescribed and realized mathematical organisation.

Our epistemological reference model is based on the contention that a similar divide can be described and observed concerning other key elements of secondary level calculus, namely derived functions and integrals. In fact, when we consider the following (rough) definition

$$f'(x) = \lim_{h \to 0} \frac{f(x+h) - f(x)}{h}$$

we see immediately that the definition of derivatives and the justification of the rules governing their behavior may indeed be considered as generating a local mathematical organisation which is directly derived from $MO_2$ as described above. More generally, the definition of the derivative generates a technology unifying a local mathematical organisation $MO_4$ whose most basic types of task are, for a given function, to describe what $f'(x)$ is, to determine whether it exists, and to justify the so-called "rules of differentiation". These "rules" also constitute the theory level of an "algebra of differentiation" $MO_3$ which, as before, can exist in relative independence from $MO_4$.

We should not fail to note here that important theoretical results in differential calculus – like the mean value theorem – rely not just on $MO_2$ but also on other local organisations unified by a theory on the real number system; and some of these results are indeed important to justify other basic elements of secondary level analysis (like the link between the derivative of *f* and the monotonicity of *f*). So even for the purpose of analyzing secondary level analysis, an epistemological reference model could not consider the theory of derivatives as merely derived from $MO_2$.

Another significant difference – not least for didactic transpositions – is that ultimately differentiation is an operation which, from a given function, produces another function – not just a number, as in the case of limits. This need to think of functions as objects is further accentuated in the case of differential equations, and has been extensively discussed in the literature on presumed cognitive obstacles to calculus (see e.g. Tall, 1997). It seems, however, plausible that it is also a didactic obstacle because the mathematical organisations encountered by students before $MO_3$ do usually not have practical blocks with functions as algebraic objects (i.e. objects "to be calculated with", and legitimate as "answers").

While an exhaustive model is not the main aim here, we contend that other local organisations of differential calculus – such as those based on optimization tasks or to the solution of differential equations – can also be described in terms of an algebraic local organisation (related to computational tasks) and a topological one (related to the definitions, conditions and justifications of what and how the computation is done).

Finally, the last "grand object" of secondary level analysis is the definite integral. Again there are two basic questions to be asked, given a function defined in an interval $I = [a, b]$: does the integral $\int_a^b f(x)dx$ exist, and if so, how do we find it? From the "academic mathematics" point of view, this is related to what Jablonka and Klisinska (2012) investigated as the meaning of "the fundamental theorem of calculus", in history as well as in the minds of contemporary mathematicians. With several possible variations in the formulation, this theorem provide answers to the two basic questions just mentioned, and states that:

(1) If *f* is continuous on *I*, then *f* has an antiderivative on *I*; and if *f* has an antiderivative on *I*, then *f* is integrable on *I*.
(2) If *F* is an antiderivative to *f* on *I*, then $\int_a^b f(x)dx = F(b) - F(a)$.

The said variations in the formulation of the theorem are less interesting for its meaning than how one defines $\int_a^b f(x)dx$ to begin with. In fact, many text books (both for secondary and tertiary level) use the conclusion of the theorem as a definition (i.e. they define the integral in terms of an antiderivative). Then, of course, the theorem disappears. Still, one has an excellent new local organisation $MO_5$, the algebra of integration, with rules that are, even at the theoretical level, easily justified from the rules at the theoretical level of $MO_3$. This also suffices for the needs of some of the more advanced local organisations of differential calculus, like the algebra of solving separable differential equations. In fact, this definition works well as long as one does not seek any separate meaning in the number $\int_a^b f(x)dx$ – or in "the fundamental theorem of calculus".

Of course most introductions of the integral also relate it to area. And in some contemporary textbooks, one finds a slightly different approach to defining the integral: for a positive function *f* it is defined as the *area* of the point set $\{(x,y): a < x < b, 0 < y < f(x)\}$ while assuming tacitly or informally that this area makes sense for "good functions". Clearly, this is just like defining the limit informally: the definition makes sense in an intuitive way, but it does not suffice to enable a mathematical practice block related to MO$_6$, such as deciding on the existence of the object defined or justifying the basic rules and properties satisfied by this "area" integral. This entrance to integrals does not, however, need to leave the link to the "integral by derivatives" entirely in the dark: if one accepts the definition of the "area" integral above as meaningful in itself, one may show from first principles that it is – when viewed as a function of *b* – an anti-derivative of the function *f*. This is indeed done in most contemporary Danish text books for upper secondary school and has undoubtedly been repeated by students thousands of times at the oral part of the mathematics exam.

Even in academic (or scholarly) mathematics, the integral is defined in different ways, and development of alternative approaches offers an interesting chapter in the history of analysis, as exposed by Jablonka and Klisinska. While Lebesgue integration is often considered superior for advanced purposes, a more basic approach is the one due to B. Riemann, and a didactic transposition of it to high school is explained in the next section. But with any rigorous definition of integrals, the topological counterpart to MO$_5$ appears on the scene: a local organisation MO$_6$ unified by the theoretical definition of the integral, linked to the fundamental properties of the real number space, and with the practical block being concerned with the tasks of deciding on the existence of the integral and with justifying the rules governing its calculus.

| OBJECT | Existence / "topology" | Computation / "algebra") |
|---|---|---|
| Limit of function *f* at point $a \in [-\infty, \infty]$. | **MO$_2$** <br> $T_{21}$: Does $\lim_{x \to a} f(x)$ exist? <br> $T_{22}$: Justify rules and properties → | **MO$_1$** <br> $T_{11}$: Find $\lim_{x \to a} f(x)$. <br> THEORY BLOCK |
| Derivative of function *f* | **MO$_4$** <br> $T_{41}$: Does $f'$ exist? Where? <br> $T_{42}$: Justify rules and properties → | **MO$_3$** <br> $T_{32}$: Find $f'$. <br> THEORY BLOCK |
| Integral of function *f* on interval $[a,b] \subseteq [-\infty, \infty]$. | **MO$_6$** <br> $T_{61}$: Does $\int_a^b f(x)dx$ exist? <br> $T_{62}$: Justify rules and properties → | **MO$_5$** <br> $T_{51}$: Find $\int_a^b f(x)dx$. <br> THEORY BLOCK |

**Table 1. A model for secondary level analysis: local organisations and basic task types.**

With this, we have extended the epistemological reference model from Barbé et al. to cover the elements of secondary level analysis (or, in some countries, the introduction to university level calculus); the result is illustrated in Table 1. The point is that introductory analysis can be roughly modeled as pairs of local mathematical organisations – algebraic and topological ones – teaming up in regional ones which build on each other more or less in the sequence shown. The algebraic organistions exhibit practical blocks with algorithmic techniques which can be taught and learned if not with ease, then at least in an orderly fashion, task type by task type (it is this part which is called "calculus" in the American text books). On the other hand, the meaning of it all is related to "topological" definitions and properties which, indeed, are also needed for a deeper justification of the "calculus", but which is less evident to transpose to the classroom because of the ultimate reliance on a complete theory of the real numbers.

We have already pointed out that the six local organisations presented above and in Figure 1 do not exhaust even the most modest version of analysis at secondary level. Also the "task types" in the table are, in reality, declined into smaller collections of tasks, each characterized by one technique. So the role of the model presented here is not to be comprehensive or give all details, but instead to help us articulate principal and crucial challenges for any didactic transposition of analysis, and in particular to support our reflections on the meaning and character of the two recent "major changes" mentioned in the introduction.

**The case of integration : A didactic transposition from the past**
The most eye catching changes of the 1961-reform of Danish high school was the introduction of elements of logics, set theory and abstract algebra. Some of these elements can be made useful also to define and study functions in the domain of analysis. As was mentioned in the introduction, the reform did not dramatically affect the tasks related to analysis which appear in the final written examinations of Danish high school, although an increase in variation and difficulty of exam tasks is evident. The novelties in abstract algebra are more visible, even in the analysis tasks (with a siginificant change in terminology from "curves" to "functions" as the objects to be examined). In terms of our epistemological reference model, the exam tasks all relate to the practice blocks of the algebraic organisations $MO_1$, $MO_3$ and $MO_5$. A typical exam exercise is the following from 1971 (Petersen og Wagner, 2003, p. 256):

> A function $f$ is given by $f(x) = xe^{-2x}$, $x \in R$, where R designates the set of real numbers.
>
> Investigate $f$ as regards its zeros, sign and monotonicity.
>
> Determine the area of the point set given by $\{(x,y) \mid 0 \le x \le \frac{1}{2} \wedge 0 \le y \le f(x)\}$
>
> For any positive real number $a$ a function $g_a$ is given by $g_a(x) = xe^{-ax}, x \in R$.
>
> Show that $g_a$ has a maximal value and find it.

As in many other tasks, the analysis appears in the "investigation" of certain properties of a given function, for instance to find its asymptotes (reduces to find one or more limits, i.e. to $MO_1$-tasks), to determine monotonicity or suprema (the key to which is to find $f'$, i.e. an $MO_3$-task), or to determine the area or volume of certain figures (reduces to a definite integral, i.e. an MO$_5$-task). All of these tasks continue to be common at the written examinations.

Text books from the period 1961-1980 reveal more profound additions to the "theory blocks" taught, and certainly also marked differences with contemporary teaching at this level. In fact, all of the six local organisations described above are covered in detail, both in exposition of theory, in worked examples and in exercises. Today's university students of mathematics usually refuse to believe that this could be done at the secondary level because it is now part of their first year. To show that and how it was really done, I will provide a rather extensive exposition of the presentation of integration in a text book series authored by Kristensen and Rindung (1973), which vastly dominated Danish high school from the late sixties to the early eighties. Here we study only the first of the two books written for the second year of high school, and only the second edition from 1973. This edition differs from the 1963 edition in several respects; most notably it has a less rigorous treatment of the topology of the real numbers. For instance, in the 1973 edition, all mention of supremum and infimum was dropped. This clearly affected also the introduction of the definite integral (in fact, the Riemann integral) which we now present.

The chapter on integration has 12 main sections (we provide a short description in parentheses):

- *Area* (8 pp., an informal discussion of area of non-polygonal point sets, and how it may be approached through double approximation with polygonal point sets)
- *Mean sums, upper sums, lower sums* (5 pp., a rigorous definition of these notions for bounded functions on an interval, ending with the theorem that every lower sum is less than every upper sum)
- *Integrability* (3 pp., rigorous definition by *the existence of a unique number situated between all lower sums and all upper sums*; proof than every monotonous function is integrable)
- *The integral and mean sums* (3 pp., proof that the integral, if it exists, is a limit of mean sums and can be considered as a "mean value" of the function on the interval)
- *Interval additivity theorem* (4pp., rigorous proof given based on the above definition)
- *The class of integrable functions* (3 pp., applies additivity theorem to prove that piecewise monotonous functions are integrable. A discussion of examples and more general results, including the theorem that continuous functions are integrable – stated without proof).
- *Integral and antiderivative* (4pp., rigorous proof that if a function is integrable and has an antiderivative, then the formula above applies).
- *Existence of antiderivatives* (2 pp.). Proof that if a function *f* is continuous on an interval *I* and $a \in I$ then $F(x) = \int_a^x f(x)dx$ is an antiderivative to *f* on *I*.
- *The indefinite integral* (1 p., introduction of the symbol $\int f(x)dx$ for the class of antiderivatives).
- *Calculation rules for integration* (15 pp., including substitution and parts, with many examples).
- *Existence of logarithm functions* (2 pp., continuation of a "gap" left in the first year volume, filling it by proving that the integral of 1/*x* gives a function with the previously stated properties).

- *Application of integral calculus* (10 pp., including volumes, curve length and examples from physics and financial theory).

As this outline shows, the text presents techniques and theory covering most of $MO_6$, with the single exception that integrability is only shown for piecewise monotonous functions, not for general continuous functions. It essentially presents this *before* $MO_5$ and approximately the same space is allowed for each of these local organisations, the main link being the justification that integrals of "common functions" can be computed by antiderivatives. For practical purposes, integrability is certainly sufficiently covered as all functions normally considered in high school are piecewise monotonous, even if the book does present one continuous functions that is not (ibid., p. 154). The need to state the theorem about integrability of continuous functions is due to its use to prove the existence of antiderivatives. This is an important point in view of $MO_5$, since there are simple and common functions to which none of the "calculations rule" succeed in producing the antiderivative – thus, unlike differentiation and limits, it appears harder to dismiss the existence problem with the notion that "we only consider functions where the algebraic rules apply".

Clearly, the text book exposition of theory from $MO_6$ does not in itself guarantee that students will engage in any related practice, besides absorbing and reciting proofs on demand. So a really interesting feature of the chapters dealing with "the topology of integrals" are the attempts to engage students in solving tasks. Here are some examples of exercises from this same text book series:

**301.** Show that $f$ given by $f(x) = \begin{cases} x, & x \text{ rational} \\ 2x, & x \text{ irrational} \end{cases}$ is not integrable on $[0,1]$.

**308.** Assume that $f$ is a bounded and integrable function on the interval $[a, b]$. The function $F$ is defined by $F(x) = \int_a^x f(t)dt$ for any $x$ in $[a, b]$. Show that $F$ is continuous.

Indeed, to solve the first exercises, students must show that every lower sum is smaller than ½, while every upper sum is at least 1 – that is, they will mobilize a genuine $MO_6$-technique (to show non-integrability based on the definition mentioned above). Similarly, the second exercise requires putting the interval additivity theorem to use, together with techniques related to inequalities (a central part of $MO_6$). It is indeed possible that both the practice and theory related to $MO_6$ was not studied with the same intensity by all classes at the time. In fact, the national written exam concerned exclusively $MO_5$ – in particular, not one exam task ever asked for the integrability of a function – and for the final oral exam, more concerned with theorems and proofs, the teacher always had some freedom to select emphases and topics. However past syllabi (e.g. Petersen and Vagner, 2003, p. 243) as well as the authors' personal memory confirm that both the theory and practice of $MO_6$ were certainly developed according to the ambitions of the text book and its task inventory. But as noted by experienced teachers (ibid., p. 266),

> over the 1970's, the students increasingly had difficulties to appreciate the cautious and stringent fashion in which topics were treated in "Kristensen and Rindung". The reason was, among other things, the "learning by doing" pedagogy that grew in importance in primary and lower secondary school…

It is evident that the "new math" period ended in a more quick and abrupt way in Danish primary and lower secondary school. In high school, the use of "Kristensen and Rindung" continued well into the 1980's; the author of this chapter remembers working the two exercises quoted above in 1984. This difference is, I believe, not unrelated to the fact that Danish teachers at the primary and lower secondary level do not study mathematics at universities and as a result, have little or no experience with modern mathematics. But there is no doubt that different external constraints on the two types of school institutions were also important.

**The case of integration : An example of recent didactic transposition**
After major reforms on the 1980's and 1990's, Danish high school has become more diversified with several streams and options, which makes it more difficult to describe a typical approach to a sector like integration. The general tendency, already alluded to above, is clearly that $MO_6$ is not taught except for mentioning the link between the definite integral and certain "areas" which are assumed to make sense as a piece of nature. Clearly, $MO_5$ has become even more dominant but it has also changed, as the use of symbolic calculators for integration is now both allowed and taught along with non-instrumented techniques that are still required in parts of the final written exam (see Drijvers, 2009 for a more detailed study of CAS-use in final high school exams in the Danish and other contexts). However, it may still be possible to study informal techniques related to $MO_6$ as so-called "optional" topics. We have chosen to present some ideas from a text book by Bregendal, Schmidt and Vestergaard (2007) which illustrate how this can be done in continuation of the "mandatory" material.

The book has two chapters on integration, the first covering the mandatory material and the second dealing with more advanced options. The latter include, nowadays, $MO_5$-techniques like integration by parts, but the chapter also features a 10 page section on "Numeric integration" which is the excerpt we will consider here, as it is the part of the book which comes closest to ask the question about the existence of integrals.

The section in question opens with an informal description of how Archimedes approximated $\pi$, the area of the unit circle, by computing the areas of inscribed and circumscribed regular polygons.

The authors go on to explain how "left and right sums", corresponding to the areas of certain rectangles, can be used to do something similar in order to compute the area of $\{(x, y) : 0 < x < 1,\ 0 < y < x^2\}$, and that in fact the average of these two sums (corresponding to trapezes through the middle points of the rectangles) gives a good approximation already for just four intervals. The "good" value ($\frac{1}{3}$) is known because it has presumably been established in the basic chapter that the area can be computed using the antiderivative. The authors then explain in great detail that the *n*'th right sums are

$$H_n = \sum_{i=1}^{n} \left(\frac{i}{n}\right)^2 \frac{1}{n} = \frac{2n^3 + 3n^2 + n}{6n^3}$$

and a similar formula for the left sums is given (to be proved in an exercise). Both converge to $\frac{1}{3}$. They then state (p.80) – with no justification – that:

> A similar relation between these sums and the integral $\int_a^b f(x)dx$ is valid for any continuous function $f$ on an interval $[a, b]$, and we can in particular conclude that for such a function, one has
>
> $$H_n = \sum_{i=1}^{n} f(x_i) \cdot \Delta x \to \int_a^b f(x)dx \text{ as } n \to \infty \text{ (that is, when } \Delta x \to 0\text{).}$$

The authors proceed to show a graph of a non-specific concave function (see Figure 1) and point out that the "comparative size of the left and right sums depends on the form of the graph" (p. 80). In fact, an attentive reading of the figure shows that neither is clearly smaller or greater than the integral. This certainly blurs the connection to the idea of Archimedes. Still, an informal connection between the limit of a sum and the area has been established for a concrete increasing function where right and left sums do enclose the "area" to be computed.

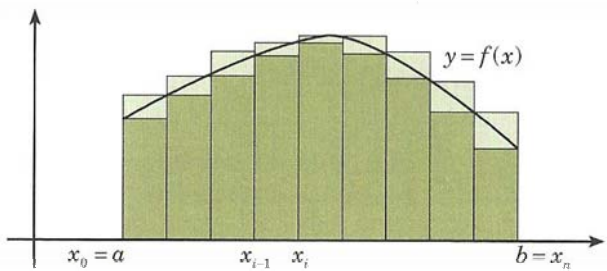

**Figur 1 : illustration of right and left sums (Bregendahl et al., 2007, p. 80)**

At the end of the section, some concrete worked examples and exercises are given on how to compute right and left sums using the statistics package of a calculator (Texas TI-84) for large values of $n$. No use of graphical visualization is suggested at this point.

In terms of our reference model, one could first think that the authors really seek to give an informal treatment of the Riemann integral; this impression is confirmed by a historical note in the book margin (p. 81), with a picture of B. Riemann and a text claiming among other things that

> The German Riemann clarified the properties of a function that make it integrable. For this reason the integral we have worked with is also called the Riemann integral.

However, nowhere else is the notion of "integrable" mentioned in the text. It is never said that the existence of the limit is a condition for the integral to exist; it is merely postulated that the formula is "valid" (quote given above). We do not get near to a distinction between lower, upper, middle, right and left sums, which would certainly be needed for a more formal treatment of the Riemann integral. Only the two last "sums" are dealt with, but they do not really correspond to the idea of Archimedes which introduced the section, or to the main idea behind Riemann's integral.

In fact, the main point seems to be to give an alternative (and potentially instrumented) technique for "computing the integral", namely the informally topological "formula" for the integral (as the limit of sums). While the algebraic techniques used to compute $H_n$ for the case $f(x) = x^2$ on $[0,1]$ will certainly not go much beyond that example, the "numeric" technique establishes a kind of experimental relation between a limit process (infinite sum) and the integral defined (and computed) using antiderivatives. The authors noted some pages earlier (p. 74) that it is sometimes impossible to find an antiderivative "using the methods we have seen", and this then motivates the introduction of "numerical methods". As the integral is originally *defined* using antiderivatives, and then "shown" to correspond to an area in some cases, the limit formula is in fact introduced as an alternative technique for *finding* integrals (i.e. in $MO_5$), not as a tool to *define* integrability and integrals (i.e. as a technique for $MO_6$).

In short, the epistemic value of the "topological technique", which is at the root of Riemanns integral (and of $MO_6$), fails to appear in the text. And the pragmatic value of the alternative numerical technique (whether instrumented or not) may be equally unconvincing to students in possession of calculators who do numerical computation of definite integrals in one step. Of course, the same can be said of most of the techniques of $MO_5$, as any symbolic integration that manual techniques can achieve, is also done in one step by the students' CAS devices.

**The dilemma – and a challenge**

Calculators became mandatory tools in Danish high school mathematics around 1980. At first, these handheld devices replaced tables and other tools to compute values of special functions etc., thus suppressing previously important techniques and tools. During the eighties, a rapid succession of more advanced calculators appeared – programmable, graphical, and "computer algebra system" calculators. With more or less delay, the use of some or all of these devices – as well as similar laptop software – has become part of high school teaching of mathematics. This is not the place to go into the historical or didactical subtleties of this development (we refer to Hoyles and Lagrange, 2010, for an excellent entry). We just stress that the CAS systems which are used in present-day high school teaching, at least in Denmark, provide instrumented techniques of high pragmatic value (or efficiency) for all basic tasks in $MO_1$, $MO_3$ and $MO_5$ cited above. Most of the customary "standard tasks" related to investigating given function are simplified, if not trivialized, using these techniques. At the same time, the connection between local organisations – such as the key connection between $MO_3$ and $MO_5$ as "opposite tasks" – tends to disappear when considering these organisations with instrumented techniques. Of course, it is possible and necessary, then, to develop new tasks for students, both for the daily teaching and final examinations. Indeed the interpretation of more or (often) less authentic situations in terms of function "models" is a much treasured direction for doing so, at least in Denmark (cf. also Drijvers, 2009). But a certain dissatisfaction in terms of the mathematical coherence cannot be denied. When the computation of limits, derivatives, maximal and minimal values, antiderivatives and so on reduced to independent, one-key operations, not much is left of the algebraic organisations and their theoretical coherence. At the same time, the topological organisations $MO_2$, $MO_4$ and $MO_6$ have been long abandoned, at least in the formally demanding transpositions they used to have.

The dilemma we then face is the following: what used to be the "core contents" of high school mathematics for decades – almost 80 years in Denmark – seems now reduced to a collection of independent, highly instrumented techniques together with a basic algebraic technology of functions and numbers, which ebables them to be used and combined to solve a variety of variation problems of real importance in many settings. The theory of computation – corresponding to rules of the theory blocks of $MO_1$, $MO_3$ and $MO_5$ – continues to be taught and learned in a more or less complete and abstract form, but their practical value is to a large extent gone at least for beginners. We have already explained the incoherence resulting from the elimination of the topological parts of mathematical analysis in the transposition to high school mathematics; with the introduction of instrumented techniques, we may face a more or less complete collapse also when it comes to the coherence which remained among and inside the local algebraic organisations.

Of course the dilemma can also be considered as a challenge: how can we reorganize – or modernize – the transposition of mathematical analysis to and in high school teaching in ways that make use of the affordances of technology while presenting the mathematical domain of analysis in a more complete and satisfactory way than as a set of modeling tools?

The answer of the "past transposition" by Kristensen and Rindung (1973) was essentially to keep as closely as possible to the "scientific" mathematics of its time; clearly, the pedagogical and political trends make that principle less evident today. But it should at least be noted that the proximity principle of the past might not lead to exactly the same answers as it did the 1960s. One reason is that CAS-based experimental methods have become part of the "scientific" practice also to the scientist who develops and uses mathematical analysis. The heart of analysis – which remains limits and deep properties of the real number systems – can be accessed and treated in new ways using technology, beginning with somewhat ostensive approaches (for instance, as applets "showing" definitions, e.g. www.maplesoft.com/products/mapleplayer/).  Another reason is that the mathematical analyst of today is yet another generation from the time where rigorous analysis was something new and exciting – functions, limits, and the other key objects have somehow been "tamed" by the mathematical practice, just like complex and negative numbers a little earlier. This could lead to a higher tolerance for relative informal approaches in the secondary curriculum, as long as the transposition preserve what the contemporary scientist regards as essential to the transposed mathematical organisations.

The answer of the "recent transposition" of the Riemann integral, presented above, appears clearly unsatisfactory, even if it hints at the potential interest of a sequence approach to the topological side of elementary analysis. The dilemma identified above is, in a way, only accentuated by adding another more or less unjustified technique to the transposition of $MO_5$. An interesting alternative would be to introduce the integral *first* as the limit of (say) right sums, with convergence being the condition of existence; *then* derive some of the properties that allow for (some sort of) proof that the derivative of the integral is the integrand. A related but much more radical alternative, is to revert the common transposition and present integration (including *both* $MO_5$ and $MO_6$) before differentiation. This approach was already completely developed by Apostol (1967) who advocated the choice by appealing to the historic precedence of problems related to integration (in fact, the Archimedean arguments alluded to above). I do not know of high school text books taking this approach, which still seems to appear almost offensive even to some college teachers (see e.g. Math Forum, 2009). However, the use of instrumented techniques for computation and visualization, and the many attractive uses of integration, could well mean that teaching this sector first might become an interesting option at the secondary level.

By way of conclusion, the teaching of analysis in secondary school is not only threatened in its time-honored form by the affordances of new technology. In fact, the exercise of certain algebraic techniques, as the main element in students' praxeologies in analysis, had already become critically separate from the mathematical and extra-mathematical questions that motivate their development and also from each other, in the absence of theoretical elements that could help to relate and justify them as mathematical practices. Research and development concerning the the secondary curriculum in analysis should not only focus on the algebraic side, despite the obvious interest of technology in this setting – it should also seek ways in which mathematical software and other resources can help "rebalance" the fundamental synergy between algebra and topology in this topic. This means, in particular, to give students access to its fundamental constructs – limits, derivatives and integrals – in ways that will make them useful tools to solve real questions involving infinite sums, mean values, growth rates and so on.

**References.**

Apostol, T. (1967). *Calculus I*, second ed. New York: John Wiley and Sons.

Barbé, J., Bosch, M., Espinoza, L. and Gascón, J. (2005). Didactic restrictions on teachers' practice – the case of limits of functions in Spanish high schools. *Educational Studies in Mathematics* 59, 235-268.

Bregendahl, P., Schmidt, S. and Vestergaard, L. (2007). *Mat A hhx*. Aarhus: Systime.

Drijvers, P. (2009). Tools and tests: technology in national final mathematics examinations. In C. Winsløw (Ed.), *Nordic Research on Mathematics Education, Proceedings from NORMA08*, pp. 225-236. Rotterdam: Sense Publishers.

Gyöngyösi, E., Solovej, J. and Winsløw, C. (2011). Using CAS based work to ease the transition from calculus to real analysis. In M Pytlak, E Swoboda & T Rowland (Eds), Proceedings of the seventh congress of the European society for research in mathematics education, Rzeszow, s. 2002-2011.

Hoyles, C. And Lagrange (Eds., 2010). *Mathematics Education and Technology-Rethinking the Terrain. The 13$^{th}$ ICMI study*. Springer: New York.

Jablonka, E. and Klisinska, A. (2012), What was and is the fundamental theorem of calculus, really? In B. Sriraman (Ed), *The Montana Mathematical Enthousiast Monographs No. 12:*

*Crossroads in the history of mathematics and mathematics education*, pp. 3-40. Missoula : Information Age Publ. & The Montana Council of Teachers of Mathematics.

Kristensen, E. And Rindung, O. (1973). *Matematik 2.1*.

MathForum (2009). http://mathforum.org/kb/message.jspa?messageID=6858531

Petersen, P. and Vagner, S. (2003). *Studentereksamensopgaver i matematik 1806-1991*. Hillerød: Matematiklærerforeningen.

Spresser, D. (1979). Placement of the first calculus course. *International Journal for Mathematics in Science and Technology* 10, 593-600.

Tall, D. (1997): Functions and calculus. In A. J. Bishop et al (Eds.), *International Handbook of Mathematics Education*, 289-325. Dordrecht: Kluwer (1997)